\documentclass[a4paper,12pt]{article}
\usepackage{amsmath,amsfonts,amssymb,latexsym,epsfig,mathrsfs}
\usepackage[all]{xy}
\usepackage{graphicx}
\newcommand{\be}{\begin{equation}}
\newcommand{\ee}{\end{equation}}

\newtheorem{example}{Example}[subsection]
\newtheorem{theorem}{Theorem}[subsection]
\newtheorem{definition}{Definition}[subsection]
\newtheorem{proposition}{Proposition}[subsection]

\newtheorem{remark}{Remark}[subsection]
\def\proof{\bf Proof. \rm}

\begin{document}
\thispagestyle{empty}
 \begin{center}
{{\bf ON THE CASE WHERE ADJOINT AND COADJOINT ORBIT SPACES ARE SYMPLECTOMORPHIC}
	
	\bigskip
	
	{\small Augustin T. Batubenge\footnote{Corresponding author} \quad \quad \quad \quad  Wallace M. Haziyu }}

\end{center}

\vspace{0.5cm}

\begin{abstract}
Let $G$ be an n-dimensional semisimple, compact and connected Lie group acting on both the Lie algebra $\mathfrak{g}$ of $G$ and its dual $\mathfrak{g}^*$. In this work it is shown that a nondegenerate Killing form of $G$ induces an $Ad^{*}$-equivariant isomorphism of $\mathfrak{g}$ onto $\mathfrak{g}^*$ which, in turn, induces by passage to quotients a symplectic diffeomorphism between adjoint and coadjoint orbit spaces of $G$.
\end{abstract}


\section{Introduction}
This work is concerned with morphisms of the category of symplectic spaces, so-called symplectic mappings. Of more interest among them are isomorphisms, that is, the symplectic mappings which also are diffeomorphisms between objects. They are important in that they exchange both the differentiable as well as symplectic structures. Working in this area so-called symplectic geometry is fascinating in that several studies, going back to previous centuries, constantly aimed at working out an elegant formalism of classical mechanics. For this paper, our main references among others are the book by R. Abraham and J.E. Marsden (\cite{Abr78}) and A. Arvanitoyeorgos (\cite{Arv03}). Combining the information provided in these sources as well as the constructions in the recent author's paper (see \cite{BH18}), we were able to obtain the main results of this study. The work involves an important amount of background ideas on representation theory. That is, adjoint and coadjoint representations as well as the actions of Lie groups yielding orbit spaces. With high interest are those quotient spaces resulting from transitive actions, the homogeneous spaces. These are the Lie groups themselves, spheres in real as well as complex and quaternionic settings, projective spaces, Grassmann and Stiefel manifolds to cite a few. In the list, we would mention flag and generalized flag manifolds. They are an important class of homogeneous spaces which admit a complex structure, a K\"ahler structure and a symplectic structure as mentioned in (\cite{Arv03}). \\The study of coadjoint orbits was introduced by Kirillov, and the existence of a symplectic structure on these orbits is the result of Kostant and Souriau (see \cite[p.52]{Ber01}), the fact from which we shall take steps further in this study. Briefly speaking, we consider the action of a compact, connected and semisimple Lie group $G$ on both its Lie algebra $\mathfrak{g}$ and its dual $\mathfrak{g}^*$, resulting into orbit spaces which consist of only one orbit each, and constructed a symplectic diffeomorphism between them. To this end, the paper is organized as follows. We begin by recalling the basics on homogeneous spaces. Then the notion of an adjoint orbit will follow, and we will show that it is related to flag as well as symplectic manifolds. Next, using Cartan's criterion for semisimplicity in which case the Killing form is nondegenerate and $Ad$-invariant, we will construct an $Ad^*$-equivariant isomorphism of Lie algebras $\mathfrak{g}$ onto $\mathfrak{g}^*$ that will induce a symplectomorphism on the quotient spaces reduced to one orbit each.        

\section{Preliminaries}

Let $G$ be a Lie group, $H$ a subgroup, and $G/H = \lbrace aH:a\in G\rbrace$ the set of left cosets of $H$ in $G$. The map $\pi:G\rightarrow G/H$ which takes each element $a\in G$ to its coset $aH$, is called the projection map. The coset space $G/H$ is not necessarily a manifold. However, if $H$ is a closed subgroup of $G$,  a manifold structure on the quotient space $G/H$ can be defined such that the projection map $\pi :G\rightarrow G/H$ is a surjective submersion. (see \cite [Theorem 9.2]{Boo75}). Also, recall that if $\phi:G\times M\longrightarrow M;~\phi(g,p)=\phi_g(p)$ is a smooth and transitive action of $G$ on a smooth manifold $M$, the $M$ is called a homogeneous space (see \cite [150]{Boo75}). This definition extends to the quotient space $G/H$ of the Lie group $G$ by a closed subgroup $H$ of $G$. In effect, there is a natural action $G\times G/H\rightarrow G/H$, $(g,aH)\mapsto gaH$. This action is always transitive since if $aH, bH\in G/H$, then $ba^{-1}(aH) = bH$ for all $a,b\in G$. For this reason every transitive action can be represented as a coset space $G/H$ where $H$ is a closed subgroup of $G$. In fact if $M$ is a manifold on which a Lie group $G$ acts transitively, then for any $p\in M$ let $G_{p} = \lbrace g\in G: g\cdot p = p\rbrace$ be the stabilizer of $p$, we have that $G_{p}$ is a closed subgroup of $G$ and $G/G_{p}\cong M$. Take $H = G_{p}$. Then $M\cong G/H$ as asserted. Therefore, $G/H$ is called the homogeneous space of $M$.\\

\section{Adjoint Orbits}

\begin{definition}
Let $G$ be a Lie group and $\mathfrak{g}\cong T_{e}G$ be its Lie algebra where $e$ is the identity element in $G$. Then the smooth action
\[
\Phi:G\times\mathfrak{g}\rightarrow\mathfrak{g};\quad (g,\xi)\mapsto Ad(g)\xi
\]
denoted by $Ad$, is called the adjoint action of $G$ on its Lie algebra $\mathfrak{g}$.
\end{definition}

\begin{definition}
Let $Ad:G\times\mathfrak{g}\rightarrow\mathfrak{g}$ be the adjoint action of a Lie group $G$ on its Lie algebra $\mathfrak{g}$ and let $\xi\in\mathfrak{g}$. We define the adjoint orbit of $\xi$ to be
\[
\begin{array}{ccc}
O_{\xi} = \lbrace Ad(g)\xi:g\in G\rbrace\subset\mathfrak{g}
\end{array}
\] 
\end{definition}
That is, if $\eta\in O_{\xi}$ then there is some $g\in G$ such that $\eta = Ad(g)\xi$. The stability group also called the isotropy group of $\xi$ is given by
\[
\begin{array}{ccc}
G_{\xi} = \lbrace g\in G: Ad(g)\xi = \xi\rbrace.
\end{array}
\]
This is a closed subgroup of $G$ (see \cite [p~16]{Cro18}). In what follows, we show that adjoint orbits can be represented as homogeneous spaces. For a similar construction (see \cite [pp~127-129]{BH18}). Define a map $\rho:O_{\xi}\rightarrow G/G_{\xi}$ by $\rho(\eta) = gG_{\xi}$ for $\eta\in O_{\xi}$ and $g\in G$ such that $\eta = Ad(g)\xi$. The map $\rho$ is well defined since if also $\rho(\eta) = hG_{\xi}$ for some $h\in G$ then $Ad(g)\xi = Ad(h)\xi$ which implies that $Ad(h^{-1})\circ Ad(g)\xi = \xi$. This gives $h^{-1}g\in G_{\xi}$ and $gG_{\xi} = hG_{\xi}$. The map $\rho$ is injective. For, let $\eta = Ad(g)\xi$, $\mu = Ad(h)\xi$ and suppose that $gG_{\xi} = hG_{\xi}$. Then $h^{-1}g\in G_{\xi}$ so that $Ad(h^{-1}g)\xi = Ad(h^{-1})\circ Ad(g)\xi = \xi$. This implies then that $\eta = Ad(g)\xi = Ad(h)\xi = \mu$. Clearly $\rho$ is surjective since for $g\in G$ and $\eta = Ad(g)\xi\in O_{\xi}$ gives $\rho(\eta) = gG_{\xi}$ by construction. If $\eta = Ad(h)\xi$ for some $h\in G$, then $G_{\eta} = Ad(h)G_{\xi}Ad(h^{-1})$. Thus, for all $g\in G$ we have
 \[
\begin{array}{ccc}
G/G_{\xi}\cong G/G_{Ad(g)\xi}
\end{array}
\]
induced by the map $g\mapsto hgh^{-1}$, which shows that the definition of $G/G_{\xi}$ does not depend on the choice of the element $\xi$ in its adjoint orbit. Thus, $G/G_{\xi}\cong O_{\xi}$. Now let  $G/G_{\xi}\cong O_{\xi}$. From the argument above, $G$ acts transitively on $G/G_{\xi}\cong O_{\xi}$ which makes it into a homogeneous space. \\

Next, let $X\in \mathfrak{g}$. Note that the vector field on $\mathfrak{g}$ corresponding to $X$, called the fundamental vector field or the infinitesimal generator of the action, is defined by
\[
\begin{array}{ccc}
X_{\mathfrak{g}}(\xi) = \frac{d}{dt}(Ad(\exp{tX})\xi)\mid_{t=0}
\end{array}
\]
We compute the tangent space to the adjoint orbit $O_{\xi}$ at $\xi$ as follows. Let $X\in\mathfrak{g}$. Let $x(t) = \exp{tX}$ be the curve in $G$ which is tangent to $X$ at $t = 0$. Then $\xi(t) = Ad(\exp{tX})\xi$ is the curve on $O_{\xi}$ such that $\xi(0) = \xi$. Let $Y\in \mathfrak{g}$, then $\langle\xi(t) , Y\rangle = \langle Ad(\exp{tX})\xi , Y\rangle$, where $\langle\cdot , \cdot\rangle$ is the natural pairing. Differentiating with respect to $t$ at $t = 0$ we get
\[
\begin{array}{cll}
\langle\xi'(0) , Y\rangle &=& \frac{d}{dt}\langle Ad(\exp{tX})\xi , Y\rangle\mid_{t=0} \\
						&=& \langle \frac{d}{dt}(Ad(\exp{tX})\xi)\mid_{t=0} , Y\rangle = \langle ad(X)\xi , Y\rangle.
\end{array}
\]
Thus $\xi'(0) = ad(X)\xi$. Therefore, the tangent space to the orbit $O_{\xi}$ at $\xi$ is given by

\[
\begin{array}{ccc}
T_{\xi}O_{\xi} = \lbrace ad(X)\xi : X\in\mathfrak{g}\rbrace
\end{array}
\]

\subsection{Adjoint orbits as flag manifolds}
The examples of adjoint orbits that will be of interest in this work are the generalized flag manifolds. These orbits are known to hold a symplectic structure. Generalized flag manifolds are homogeneous spaces which can be expressed in the form $G/C(S)$, where $G$ is a compact Lie group and\\
$C(S)=\lbrace g\in G : gx = xg,~ \textrm{for~all~}x\in S\rbrace$ is the centraliser of a torus $S$ in $G$. Generalized flag manifolds just like flag manifolds are homogeneous spaces (see\cite [p~70]{Arv03}). Here is an example in $\mathbb{C}^n$.

\begin{definition}
Let $\mathbb{C}^{n}$ be an $n$-dimensional complex space. A flag is a sequence of complex subspaces 
\[
\begin{array}{cc}
W = V_{1}\subset V_{2}\subset\cdots\subset V_{n} = \mathbb{C}^{n}
\end{array}
\]
ordered by inclusion such that $\dim V_{i} = i$ for $i = 1,\cdots,n$ and $V_i$ is a proper subset of $V_{i+1}$ for $i=1,..., n-1.$
\end{definition}
\begin{example}
Let $\lbrace e_{1}, e_{2},\cdots , e_{n}\rbrace$ be the canonical basis for the complex vector space $\mathbb{C}^{n}$. Then the standard flag is given by
\[
\begin{array}{cc}
W_{0} = Span_{\mathbb{C}}\lbrace e_{1}\rbrace \subset Span_{\mathbb{C}}\lbrace e_{1} , e_{2}\rbrace\subset\cdots\subset Span_{\mathbb{C}}\lbrace e_{1},\cdots e_{n}\rbrace = \mathbb{C}^{n}
\end{array}
\]
\end{example}

We need to show that flag manifolds are homogeneous spaces. Let $F_{n}$ be the set of all flags in $\mathbb{C}^{n}$ and let $W_{0}$ be the standard flag above. Then the action of the Lie group $U(n)=\{A\in GL(n,\mathbb{C}):\bar{A}^TA=I\}$ on $F_{n}$ is transitive. For, consider an arbitrary flag
$W = V_{1}\subset V_{2}\subset\cdots\subset V_{n} = \mathbb{C}^{n}$. Then $U(n)$ acts on $F_{n}$ by left multiplication. That is, if $S\in U(n)$ then $SW=SV_{1}\subset SV_{2}\subset\cdots\subset SV_{n} = \mathbb{C}^{n}$. Start with $v_{1}$, a unit vector in $V_{1}$ such that $V_{1} = Span_{\mathbb{C}}\lbrace v_{1}\rbrace$. Next choose a unit vector $v_{2}$ in $V_{2}$ orthogonal to $V_{1}$ such that $V_{2} = Span_{\mathbb{C}}\lbrace v_{1},v_{2}\rbrace$. Having chosen unit vectors $ v_{1},\cdots, v_{k}$ with $V_{k}=Span_{\mathbb{C}}\lbrace v_{1},\cdots,v_{k}\rbrace$, choose further a unit vector $v_{k+1}$ in $V_{k+1}$ orthogonal to $V_{k}$ such that $V_{k+1} = Span_{\mathbb{C}}\lbrace v_{1},\cdots, v_{k+1}\rbrace$. Continuing this construction we obtain a set of orthonomal unit vectors $\lbrace v_{1}, \cdots, v_{n-1}\rbrace$ such that $V_{j} = Span_{\mathbb{C}}\lbrace v_{1},\cdots, v_{j}\rbrace$. Let $v_{n}$ be a unit vector in $V_{n}$ orthogonal to $V_{n-1}$. The set $\lbrace v_{1},v_{2},\cdots, v_{n}\rbrace$ is another orthonormal basis for $\mathbb{C}^{n}$. It is now a result of linear algebra that there is $n\times n$ matrix $S=(a_{ij})$ such that $v_{i}=\displaystyle\sum_{j=1}^{n}a_{ij}e_{j}$. Then $S\in U(n)$ and $SW_{0} = W$. Thus $U(n)$ acts transitively on $F_{n}$ as earlier claimed.\\

The isotropy subgroup of $W$ is $ \lbrace A\in U(n): AV_{j} = V_{j}\rbrace$. In particular, this is a set of matrices $A\in U(n)$ such that $Av_{k} = \lambda_{k}v_{k}$ for some complex number $\lambda_{k}$ with $\mid \lambda_{k}\mid = 1$ since $A\in U(n)$. Thus $\lambda_{k} = e^{i\theta_{k}}\in U(1)$. Since this must be true for each $v_{j}$, $j=1, 2,\cdots, n$, the matrix $A$ must be of the form $A = diag(e^{i\theta_{1}},\cdots, e^{i\theta_{n}})$. Thus $F_{n}=U(n)/U(1)\times\cdots\times U(1)$\\

Now let $\lbrace n_{1},\cdots, n_{k}\rbrace$ be a set of positive integers such that $n_{1}+n_{2}+\cdots +n_{k}=n$. A partial flag is an element $W=V_{1}\subset \cdots\subset V_{k}$ with $\dim V_{k} = n_{1}+\cdots +n_{k}$. We can visualize this as a sum of vector spaces. For example, let $Q_{1}, Q_{2},\cdots, Q_{n}$ be a set of subspaces of $\mathbb{C}^{n}$ with $\dim Q_{1}=n_{1}$ , $\dim Q_{2} = n_{2}\cdots \dim Q_{n-1} = n-1$. \\
Set
\[ 
\begin{array}{cll}
	V_{1} &=& Q_{1} \\
	V_{2} &=& Q_{1}\oplus Q_{2} \\
    	  &\cdots& \\
	V_{n-1} &=& Q_{1}\oplus Q_{2}\oplus\cdots\oplus Q_{n-1} \\
\end{array}
\]

Then $V_{1}\subset\cdots\subset V_{n-1}$ and $\dim V_{j}=n_{1}+\cdots +n_{j}$. The flag $W=V_{1}\subset\cdots\subset V_{k}$ with $\dim V_{k}=n_{1}+\cdots +n_{k}$ is called a partial flag.\\

A generalized flag manifold in $\mathbb{C}^n$ is a set $F(n_{1},\cdots,n_{k})$ of all partial flags with $n_{1}+n_{2}+\cdots +n_{k} = n$. Throughout the discussion that follows, the Lie group $G$ will be compact and connected. We chose the unitary group $U(n)$ in order to illustrate that. (see Batubenge et.al. \cite{BBK85})

\begin{description}
\item   (i) $U(n)$ is compact. \\
This is because $U(n)$ is both closed and bounded in $GL(n,\mathbb{C})$. For, $U(n)=det^{-1}(S^1)=det(U(1)),$ where we denoted by $det$ the determinant function. Next, we show that $U(n)$ is bounded. For, pick $A=(\alpha_{ij})\in U(n)$. One has $\displaystyle\sum_{j}\alpha_{ij}\cdot \beta_{jk}=\delta_{ik}$, the Kronecker delta, with $\beta_{jk}=\bar{\alpha}_{kj}$. Hence, for $i=k$ one has $\displaystyle\sum_{j}\alpha_{ij}\cdot \bar{\alpha}_{ji}=1.$ Hence, $$\displaystyle\sum_{i=1}^{n}\bigg(\displaystyle\sum_{j=1}^{n}|\alpha_{ij}|^2\bigg)=n.$$ Now, 
$$||A||=\displaystyle\bigg(\sum_{i,j=1}^n|\alpha_{ij}|^2\bigg)^{\frac{1}{2}}=\sqrt{n}<\sqrt{n+1}.$$ Therefore, one has $A\in B(0,\sqrt{n+1})$, where $r=\sqrt{n+1}.$ Now one has that $A\in B(0,r)$ whenever $A\in U(n)$ so that $U(n)\subset B(0,r)$, with $r=\sqrt{n+1}$. Hence, $U(n)$ is bounded. Thus, $U(n)$ is compact.

\item   (ii)    $U(n)$ is connected

Consider the action of $U(n)$ on $\mathbb{C}^{n}$ given by $(A,X)\mapsto AX$ for all\\
 $A\in U(n)$ and $X\in\mathbb{C}^{n}$. We have
  $$\| AX\|^{2}=(\bar{AX}^{T})(AX)=\bar{X}^{T}\bar{A}^{T}AX = \bar{X}^{T}X = \|X\|^{2}.$$ Thus, this action takes sets of the form \\ $\lbrace(z_{1},\cdots,z_{n}):\mid z_{1}\mid^{2}+\mid z_{2}\mid^{2}+\cdots+\mid z_{n}\mid^{2} = 1\rbrace$ into sets of the same kind. In particular, the orbit of $e_{1}$ under this action is the unit sphere $S^{2n-1}$. The stabilizer of the same element $e_{1}$ are matrices of the form
\begin{displaymath}
\left(\begin{array}{ccc}
1&0\\
0&A_{1}
\end{array}\right)
\end{displaymath}

where $A_{1}\in U(n-1)$. Thus $S^{2n-1} = U(n)/U(n-1)$. But $S^{2n-1}$ is connected which implies that $U(n)$ is connected if and only if $U(n-1)$ is connected. Since $U(1) = S^{1}$ is connected, we conclude by induction on $n$ that $U(n)$ is connected.
\end{description}

The Lie algebra of $U(n)$ is the space of all skew-Hermitian matrices\\
 $\mathfrak{u}(n)=\lbrace A\in Mat_{n\times n}(\mathbb{C}): A+\bar{A}^{T} = 0\rbrace$. We now want to determine the orbits of adjoint representation of the Lie group $G = U(n)$ on its Lie algebra $\mathfrak{g} = \mathfrak{u}(n)$.\\

Let $Ad:G\times\mathfrak{g}\rightarrow\mathfrak{g}$ be the action of $G$ on its Lie algebra $\mathfrak{g}$. Let $X\in\mathfrak{g}$, then the orbit of $X$ is given by\\
\[
\begin{array}{cll}
O_{X} &=&\lbrace Ad_{g}X:g\in G\rbrace\\

    &=& \lbrace Y\in\mathfrak{g}:Y = gXg^{-1} ~\textrm{for~some~} g\in G\rbrace
\end{array}
\]

This is a set of similar matrices since the action is by conjugation. Recall that every skew Hermitian matrix is diagonalizable and that all the eigenvalues of a skew Hermitian matrix are purely imaginary. This means that $X$ is $U(n)-$ conjugate to a matrix of the form $X_{\lambda} = Diag(i\lambda_{1},i\lambda_{2},\cdots,i\lambda_{n})$ for $\lambda_{j}\in\mathbb{R},\hspace{0.4cm} j=1,\cdots,n$. Since similar matrices have same eigenvalues, without loss of generality we can describe the adjoint orbit of $X$ to be the set of all skew Hermitian matrices with eigenvalues $i\lambda_{1},i\lambda_{2},\cdots,i\lambda_{n}$. Denote this set of eigenvalues by $\lambda$ and the orbit determined by the corresponding eigenspaces by $H(\lambda)$. Note that $H(\lambda)$ is a vector space since it is a closed subgroup of a linear group $GL(n,\mathbb{C})$.\\

Case 1 :    All the $n$ eigenvalues are distinct\\

Let $x_{j}$ be the eigenvector corresponding to the eigenvalue $i\lambda_{j}$, then we have $gx_{j} = i\lambda_{j}x_{j}$. This gives a 1-dimensional subspace $P_{j}$ of $\mathbb{C}^{n}$ which is a line in the complex plane passing through the origin.\\
Assuming $\lambda_{1}<\lambda_{2}<\cdots <\lambda_{n}$. Note that the eigenvectors corresponding to distinct eigenvalues are orthogonal. Now each element in $H(\lambda)$ has same eigenvalues $i\lambda_{1},\cdots, i\lambda_{n}$, however, it is only distinguished by its corresponding eigenspaces $P_{1},\cdots, P_{n}$. Thus for each $n-$tuple $(P_{1},P_{2},\cdots, P_{n})$ of complex lines in $\mathbb{C}^{n}$ which are pairwise orthogonal, there will be an associated element $h\in H(\lambda)$ and each element $h\in H(\lambda)$ determines a family of eigenspaces $(P_{1},P_{2},\cdots, P_{n})$. \\

Let $(P_{1},\cdots, P_{n})\mapsto P_{1}\subset P_{1}\oplus P_{2}\subset\cdots\subset P_{1}\oplus P_{2}\oplus\cdots\oplus P_{n}=\mathbb{C}^{n}$ and define the vector space $V_{j}$ by $V_{j} = P_{1}\oplus\cdots\oplus P_{j}$. Then $W=V_{0}\subset V_{1}\subset\cdots\subset V_{n}=\mathbb{C}^{n}$ is a flag we have already seen and the totality of such flags $F_{n} = U(n)/U(1)\times\cdots\times U(n)$ is the flag manifold described earlier. There is a bijection from $H(\lambda)$ to $F_{n}$ which associates to each element $h\in H(\lambda)$ the subspaces $V_{j} = P_{1}\oplus\cdots\oplus P_{j}$ where $P_{j}$ is the eigenspace of $h$ corresponding to the eigenvalue $i\lambda_{j}$. This shows that the adjoint orbits are diffeomorphic to flag manifolds.\\

Case 2:  There are $k<n$ distinct eigenvalues.\\

We again order the eigenvalues $\lambda_{1}<\cdots <\lambda_{k}$. Let $n_{1}, n_{2},\cdots, n_{k}$ be their multiplicities respectively. Let $Q_{j}$ be the eigenspace corresponding to the eigenvalue $i\lambda_{j}$. We assume that $\dim Q_{i} = n_{i},\hspace{0.4cm} i=1,\cdots,k$. Then the orbit of $X$ is again determined by the eigenspaces $Q_{1},\cdots, Q_{k}$. We form an increasing sequence ordered by inclusion as before

$(Q_{1}, Q_{2},\cdots, Q_{k})\mapsto Q_{1}\subset Q_{1}\oplus Q_{2}\subset\cdots\subset Q_{1}\oplus\cdots\oplus Q_{k} = \mathbb{C}^{n}$.\\

Let $F(n_{1},n_{2},\cdots, n_{k})$ be the set of all such sequences. Then the orbit of $X$ is diffeomorphic to the homogeneous space
$F(n_{1},\cdots, n_{k})=U(n)/(U(n_{1})\times\cdots\times U(n_{k}))$ which as we have already seen is a generalized flag manifold. For the variation proof of this (see \cite [Proposition II.1.15]{Aud04}).

\begin{definition}\emph{(\bf Killing form)}\\
	Let $\mathfrak{g}$ be any Lie algebra. The Killing form of $\mathfrak{g}$ denoted by $B$, is a bilinear form $B:\mathfrak{g}\times \mathfrak{g}\longrightarrow \mathbb{R}$ given by 
	$$B(X,Y)=tr(ad(X)\circ ad(Y)),\rm{for~all~}X,Y\in \mathfrak{g}$$
	where $tr$ refers to the usual trace of a mapping. 
\end{definition}
\begin{remark}  We shall call $B$ the Killing form of the Lie group $G$ provided $\mathfrak{g}$ is the Lie algebra of the Lie group $G$, in which case the Killing form $B$ is $Ad$-invariant. That is, 
	$$B(X,Y)=B(Ad(g)X,Ad(g)Y ) $$
	for all $g\in \mathfrak{g}$. (see \cite[proposition 2.10]{Arv03}). 
\end{remark}	

We further recall that by Cartan's criterion for semisimplicity, a finite dimensional Lie group $G$ is said to be semisimple if its Killing form is nondegenerate (see \cite[p. 34]{Arv03}). This criterion will play a key role in the next section. We would mention that the consequences of this criterion are as follows. Let $G$ be an $n$-dimensional semisimple Lie group. If $G$ is compact then its Killing form is negative definite. Moreover, if $G$ be an $n$-dimensional connected Lie group and the Killing form of $G$ is negative definite on $\mathfrak{g}$, then $G$ is compact and semisimple.

\subsection{Adjoint orbits as symplectic manifolds}

We have seen that the adjoint orbits of flag manifolds are determined by the eigenspaces corresponding to a set of eigenvalues $i\lambda_{1},\cdots, i\lambda_{k}$. Denote this set of eigenvalues by $\lambda$ and the orbit determined by the corresponding eigenspaces by $H(\lambda)$. Let $G=U(n)$ be a Lie group and $\mathfrak{g}=\mathfrak{u}(n)$ its Lie algebra. First note that the dimension of orbit $H(\lambda)$ is $n^{2}-n$ which is even. \\

For $X\in\mathfrak{g}$ we have seen that if $x(t) = \exp{tX}$ is a curve in $G$ tangent to $X$ at $t = 0$, then $\xi(t)=Ad_{x(t)}\xi = Ad_{\exp{tX}}\xi$ is a curve in $H(\lambda)$ passing through $\xi\in\mathfrak{u}(n)$. Then the tangent vector to this curve at $t = 0$ is given by
\[
\begin{array}{ccc}
\xi'(t) = \frac{d}{dt}Ad_{\exp{tX}}\xi\mid_{t=0}~\textrm{or}~\xi'(0) = ad(X)\xi = [\xi , X]

\end{array}
\]

We shall now construct a symplectic 2-form on the orbit $H(\lambda)$. Let $h$ be an element of $\mathfrak{u}(n)$. Define a map
\[
\omega_{h}:\mathfrak{g}\times\mathfrak{g}\rightarrow\mathbb{R};\quad
\omega_{h}(X,Y) = B(h,[X,Y])
\]

where $B$ is the Killing form of $\mathfrak{g}$, the Lie algebra of $G$.

\begin{proposition}
Let $\omega_{h}$ be as defined above. Then
\begin{description}
\item   (i) $\omega_{h}$ is skew symmetric bilinear form on $\mathfrak{g}=\mathfrak{u}(n)$
\item   (ii)    $\ker\omega_{h} = \lbrace X\in\mathfrak{u}(n):[h,X]=0\rbrace$
\item   (iii)   $\omega_{h}$ is $G$-invariant. That is, for each $g\in G$ we have\\
$\omega_{Ad(g)(h)}(Ad_{g}X,Ad_{g}Y)=\omega_{h}(X,Y)$
\end{description}
\end{proposition}

\proof Part (i) follows from the properties of the Lie bracket. For part (ii) (see \cite [p~19]{Ale96}). We prove part (iii).
\[
\begin{array}{cll}
\omega_{Ad(g)(h)}(Ad_{g}X,Ad_{g}Y) &=& B(Ad_{g}h,[Ad_{g}X,Ad_{g}Y])\\
        &=& B(Ad_{g}h,[gXg^{-1},gYg^{-1}])\\
        &=& B(Ad_{g}h,\lbrace gXYg^{-1}-gYXg^{-1}\rbrace)\\
        &=& B(Ad_{g}h,g[X,Y]g^{-1})\\
        &=& B(Ad_{g}h,Ad_{g}[X,Y])\\
        &=& B(h,[X,Y])\\
        &=& \omega_{h}(X,Y)
\end{array}
\]

Now for $h\in\mathfrak{u}(n)$ we consider the orbit map
\[
\Phi_{h}:U(n)\rightarrow\mathfrak{u}(n);\quad g\mapsto ghg^{-1}
\]

That is
\[
\Phi_{h}:U(n)\rightarrow H(\lambda)\subset\mathfrak{u}(n)
\]
Then we have $ T_{I}\Phi_{h}:\mathfrak{u}(n)\rightarrow T_{h}H(\lambda)$. But the tangent space on the orbit is generated by the vector field $ad(X)\xi = [X,\xi]$, with  $X,\xi\in\mathfrak{g}$. Define a 2-form $\Omega_{h}$ on $T_{h}H(\lambda)$ by the formula
\[
\Omega_{h}([h,X],[h,Y]) = \omega_{h}(X,Y),\hspace{0.4cm} \textrm{for}~ X,Y\in\mathfrak{u}(n)
\]

\begin{proposition}
The $\Omega_{h}$ defined above is a closed and nondegenerate 2-form on the orbit $H(\lambda)$.
\end{proposition}

\proof  First note that $\Omega_{h}$ does not depend on the choice of $X,Y\in\mathfrak{u}(n)$ since if $Z\in\ker\omega_{h}$ then we have
\[
\begin{array}{cll}
\Omega_{h}([h, X+Z],[h, Y+Z]) &=& \omega_{h}(X+Z,Y+Z) \\
				    &=& B(h,[X+Z,Y+Z]) \\
 					&=& B(h,[X,Y]+[X,Z]+[Z,(Y+Z)]) \\
 					&=& B(h,[X,Y])+B(h,[X,Z])+B(h,[Z,(Y+Z)])\\
					&=& \omega_{h}(X,Y)+\omega_{h}(X,Z)+\omega_{h}(Z,(Y+Z)) \\
					&=& \omega_{h}(X,Y) \\
					&=& \Omega_{h}([h,X],[h,Y])
\end{array}
\]

Thus, $\Omega_{h}$ is well defined. It is skew-symmetric bilinear form and $G$-invariant by the construction so it is smooth. Since the Killing form $B$ is nondegenerate, $\Omega_{h}$ is nondegenerate. We only have to show that it is closed.\\

From the formula (1) in Berndt R. (see \cite [p~73] {Ber01}) we have
\[
\begin{array}{cll}
d\omega(X,Y,Z)  &=& L_{X}\omega(Y,Z)-L_{Y}\omega(X,Z)+L_{Z}\omega(X,Y)  \\
				&+& \omega(X,[Y,Z])-\omega(Y,[X,Z])+ \omega(Z,[X,Y]),
\end{array}
\]  let $X,Y,Z\in\mathfrak{u}(n)$. Then

\[
\begin{array}{cll}
	
d\Omega_{h}([h,X],[h,Y],[h,Z])  &=& d\omega_{h}(X,Y,Z) \\
								&=& \lbrace L_{X}\omega_{h}(Y,Z)-L_{Y}\omega_{h}(X,Z)+L_{Z}\omega_{h}(X,Y)\rbrace \\
								&+&\lbrace\omega_{h}(X,[Y,Z])-\omega_{h}(Y,[X,Z]) +\omega_{h}(Z,[X,Y])\rbrace 
\end{array}
\]
We now apply the Jacobi identity to each bracket given by the braces. The second bracket gives

\[
\begin{array}{cll}

 \omega_{h}(X,[Y,Z]) &-&\omega_{h}(Y,[X,Z]) + \omega_{h}(Z,[X,Y])  \\
  					&=& B(h,[X,[Y,Z]])-B(h,[Y,[X,Z]])+B(h,[Z,[X,Y]])  \\
					&=& B(h,[X,[Y,Z]]-[Y,[X,Z]]+[Z,[X,Y]])
\end{array}
\] 

and the term in the bracket is zero by the Jacobi identity since $\mathfrak{u}(n)$ is a Lie algebra. To deal with the first bracket we have

\[
\begin{array}{cll}
L_{X}\omega_{h}(Y,Z) &=& \omega_{h}(Z,[X,Y])-\omega_{h}(Y,[X,Z]) \\
L_{Y}\omega_{h}(X,Z) &=& \omega_{h}(Z,[Y,X])-\omega_{h}(X,[Y,Z]) \\
L_{Z}\omega_{h}(X,Y) &=& \omega_{h}(Y,[Z,X])-\omega_{h}(X,[Z,Y])

\end{array}\]

Substituting into the first bracket and simplifying gives

\[
\begin{array}{cll}

L_{X}\omega_{h}(Y,Z) &-&L_{Y}\omega_{h}(X,Z)+L_{Z}\omega_{h}(X,Y) \\
					&=& 2\left( \omega_{h}(X,[Y,Z])+\omega_{h}(Y,[Z,X])+\omega_{h}(Z,[X,Y])\right)

\end{array}
\] 
which again vanishes by Jacobi identity. Thus, $d\Omega_{h} = 0$ proving that $\Omega_{h}$ is indeed closed on the orbits of the adjoint action of the Lie group $G$ on its Lie algebra $\mathfrak{g}$.

\section{Coadjoint Orbits}

We now describe briefly the orbits of the coadjoint action of a Lie group $G$ on the dual of its Lie algebra. There are many references to this section such as Abraham and Marsden (\cite {Abr78}) as well as Vilasi (\cite{Vil01}).\\

Consider the Lie group $G$ acting on itself by left translation $L_{g}:G\rightarrow G$ given by $h\mapsto gh$ for $g\in G$. This map is a diffeomorphism. So, by lifting of diffeomorphisms, induces a symplectic action on its cotangent bundle
\[
\Phi:G\times T^{*}G\rightarrow T^{*}G; \quad
    (g,\alpha_{h})\mapsto \Phi(g,\alpha_{h})=L_{g^{-1}}^{*}(\alpha_{h})
\]

This action has a momentum mapping which is equivariant with the coadjoint action. The momentum mapping of this action is given by
\[
\mu:T^{*}G\rightarrow\mathfrak{g}^{*}; \quad
    \mu(\alpha_{g})\xi = \alpha_{g}(\xi_{G}(g))= \alpha_{g}(R_{g})_{*{e}}\xi = (R_{g}^{*}\alpha_{g})\xi
\]
for all $\xi\in\mathfrak{g}$, where $\mathfrak{g}^{*}$ is the dual to the Lie algebra of $G$..\\

That is, $\mu(\alpha_{g}) = R_{g}^{*}\alpha_{g}$. Every point $\beta\in\mathfrak{g}^{*}$ is a regular value of the momentum mapping $\mu$ (see \cite [p~282]{Vil01}). So we have for each $\beta\in\mathfrak{g}^{*}$
\[
\begin{array}{cll}

\mu^{-1}(\beta) &=& \lbrace \alpha_{g}\in T^{*}G: \mu(\alpha_{g}) = \beta\rbrace\\

        &=&   \lbrace \alpha_{g}\in T^{*}G: R_{g}^{*}\alpha_{g}\xi = \beta\cdot\xi~\textrm{for~all~} \xi\in\mathfrak{g}\rbrace

\end{array}
\]
In particular, $R_{e}^{*}\alpha_{e}\xi = \beta\cdot\xi$ implying that $\alpha_{e} = \beta$. Denote this 1-form by $\alpha_{\beta}$ so that
\[
\begin{array}{ccc}
\alpha_{\beta}(e) = \beta \hspace{2cm} (1)
\end{array}
\]

For $g\in G$, applying the right translation $R_{g^{-1}}^{*}$ to Equation (1) gives a right-invariant 1-form on $G$
\[
\begin{array}{ccc}
\alpha_{\beta}(g) = R_{g^{-1}}^{*}\beta \hspace{2cm} (2)
\end{array}
\]

But now for all $g\in G$ we have
\[
\begin{array}{ccc}
\mu(\alpha_{\beta}(g)) = \mu(\alpha_{g}) = R_{g}^{*}R_{g^{-1}}^{*}\beta = \beta.
\end{array}
\]

Thus, Equation (2) defines all and only points of $\mu^{-1}(\beta)$. Since the action is defined by $\Phi(g,\alpha_{h})= L_{g^{-1}}^{*}(\alpha_{h})$, the isotropy subgroup of $\beta$ is
\[
\begin{array}{ccc}
G_{\beta} = \lbrace g\in G:L_{g^{-1}}^{*}(\alpha_{\beta})= \beta\rbrace
\end{array}
\]

From the map
\[
L_{g^{-1}}^{*}: (h,\alpha_{\beta}(h)) \longrightarrow  (gh,\alpha_{\beta}(gh))
\]
we see that $G_{\beta}$ acts on $\mu^{-1}(\beta)$ by left translation on the base points. This action is proper (see \cite [p~283]{Vil01}). Since $\beta$ is also a regular value of the momentum mapping $\mu$, then $\mu^{-1}(\beta)/G_{\beta}$ is a symplectic manifold. There is a diffeomorphism

\[
\mu^{-1}(\beta)/G_{\beta}\simeq G\cdot\beta = \lbrace Ad_{g^{-1}}^{*}\beta :g\in G\rbrace\subset\mathfrak{g}^{*}~~\textrm{(see \cite [p~284]{Vil01})}
\]
of the reduced space $\mu^{-1}(\beta)/G_{\beta}$ onto the coadjoint orbit of $\beta\in\mathfrak{g}^{*}$. Thus the coadjoint orbit $G\cdot\beta$ is a symplectic manifold. The symplectic 2-form is given by the Kirillov-Kostant-Souriau form
\[
\begin{array}{ccc}
\omega_{\beta}(\nu)(\xi_{\mathfrak{g}^{*}}(\nu),\eta_{\mathfrak{g}^{*}}(\nu)) = -\nu\cdot[\xi,\eta]~ \textrm{(see \cite [pp~302-303]{Abr78})},
\end{array}
\]
where $\xi, \eta\in \mathfrak{g}$ and $\nu\in\mathfrak{g}^{*}$.

If $G$ is semisimple, it is known that in this case, $H^{1}(\mathfrak{g},\mathbb{R}) = 0$. (See \cite [p~19]{Ale96}). Thus, if $\omega$ is closed then it is exact. So, there is a 1-form $\alpha\in\mathfrak{g}^{*}$ such that $d\alpha = \omega$. The 1-form $\alpha$ satisfies $d\alpha(X,Y) = \alpha([X,Y])$.\\

Thus if the Lie group $G$ is semisimple, compact and connected, then we have the relation\\

$\alpha([X,Y])=d\alpha(X,Y)=\omega(X,Y)=B([\xi,X],Y)=B(\xi,[X,Y])$, where $\alpha\in\mathfrak{g}^{*}$, $\omega$ a 2-form on the homogeneous space $G/H$, $B$ the Killing form on $G/H$ and $\xi,X,Y\in\mathfrak{g}$, the Lie algebra of $G$.

\section{Main results}
\begin{theorem}\label{th 511:th 511}
Let $Ad:G\times\mathfrak{g}\rightarrow\mathfrak{g}$ be an adjoint action of an $n$-dimensional semisimple, compact, connected Lie group $G$ on its Lie algebra $\mathfrak{g}\cong T_{e}G$. Let $\mathfrak{g}^{*}$ be the dual of $\mathfrak{g}$. Then there is an $Ad^{*}$-equivariant isomorphism $B^{\flat}:\mathfrak{g}\rightarrow\mathfrak{g}^{*}$.
\end{theorem}

\proof  Let
\[
B^{\flat}:\mathfrak{g}\rightarrow\mathfrak{g}^{*};\quad X\mapsto B^{\flat}(X):\mathfrak{g}\rightarrow\mathbb{R},\quad Y\mapsto B^{\flat}(X)Y:= B(X,Y)
\]   
where $B$ is the Killing form. Then $B^{\flat}$ is linear since of  for all $X,Y,Z\in\mathfrak{g}$ and using the fact that the Killing form $B$ is bilinear, we have\\
\[
\begin{array}{cll}
B^{\flat}(aX+bY)Z 	&=& B(aX+bY, Z) \\
					&=& aB(X,Z)+bB(Y,Z) \\
					&=& aB^{\flat}(X)Z+bB^{\flat}(Y)Z \\
					&=& (aB^{\flat}(X)+bB^{\flat}(Y))Z.

\end{array}
\] 
Thus $B^{\flat}(aX+bY)=aB^{\flat}(X)+bB^{\flat}(Y)$.\\
First, $B^{\flat}$ is injective. For, let $B^{\flat}(X) = B^{\flat}(Y)$. Then for all $Z\in\mathfrak{g}$ one has 
$$B^{\flat}(X)Z = B^{\flat}(Y)Z\Rightarrow B(X,Z)=B(Y,Z)\Rightarrow B(X-Y,Z)=0$$ 
and since the Killing form is nondegenerate we get $X = Y$. Next, $B^{\flat}$ is surjective since, first we note that $G$ is finite dimensional Lie group and $B^{\flat}$ is injective, thus $\ker B^{\flat} = \lbrace 0\rbrace$ implying that $\dim\ker B^{\flat} = 0$. But $\dim\ker B^{\flat}+\textrm{Rank} B^{\flat} = \dim\mathfrak{g}$, so we must have $\dim\mathfrak{g}^{*} = \dim \textrm{Im}B^{\flat} = \textrm{Rank} B^{\flat} = \dim\mathfrak{g}$. This shows that the map $B^{\flat}$ is surjective.

We now show that $B^{\flat}:\mathfrak{g}\rightarrow\mathfrak{g}^{*}$ is equivariant with respect to the adjoint action of $G$ on $\mathfrak{g}$ and the coadjoint action of $G$ on $\mathfrak{g}^{*}$. Define a map
\[
 u:G\times\mathfrak{g}\rightarrow G\times\mathfrak{g}^{*};\quad (g,X)\mapsto (g,B^{\flat}X),
\]
where  $X\in\mathfrak{g}, g\in G$. That is, $u = Id_{G}\times B^{\flat}$. Then the following diagram commutes
\[
\xymatrixcolsep{4pc}\xymatrixrowsep{4pc}
\xymatrix{
G\times \mathfrak{g} \ar[d]_{Ad} \ar[r]^-{u} & G\times \mathfrak{g}^{*} \ar[d]^{Ad^{*}}\\
\mathfrak{g} \ar[r]^{B^{\flat}} & \mathfrak{g}^{*}}
\]

Let $(g,X)\in G\times \mathfrak{g}$. Then for all $Y\in\mathfrak{g}$ we have\\

$B^{\flat}(Ad_{g}X)Y = B(Ad_{g}X,Y) = B(Ad_{g^{-1}}\circ Ad_{g}X, Ad_{g^{-1}}Y)\\
= B(X, Ad_{g^{-1}}Y) = Ad_{g}^{*}B^{\flat}(X)(Y)$. The second and the third equalities is because the Killing form $B$ is Ad-invariant. That is, 
$$B^{\flat}(Ad_{g}X) = Ad_{g}^{*}B^{\flat}X.$$

Thus $B^{\flat}\circ Ad = Ad^{*}\circ B^{\flat}$ and $B^{\flat}$ is equivariant. \\

Let $\pi_{\mathfrak{g}}:\mathfrak{g}\rightarrow \mathfrak{g}/G$ and $\pi_{\mathfrak{g}^{*}}:\mathfrak{g}^{*}\rightarrow \mathfrak{g}^{*}/G$ be the projection maps into the respective orbit spaces. Then, (see \cite [p~10]{Mei03}) there is at most one manifold structure on $\mathfrak{g}/G$ respectively on $(\mathfrak{g}^{*}/G)$ such that $\pi_{\mathfrak{g}}$ respectively $(\pi_{\mathfrak{g}^{*}})$ are submersions. In fact note for example that the rank of $d\pi_{\mathfrak{g}}$ is equal to the dimension of its image and since $\dim\mathfrak{g}/G\leq\dim\mathfrak{g}$ then $\pi_{\mathfrak{g}}$ is a submersion. Since $B^{\flat}:\mathfrak{g}\rightarrow \mathfrak{g}^{*}$ is equivariant and $\pi_{\mathfrak{g}}$ and $\pi_{\mathfrak{g}^{*}}$ are submersions, the criterion of passage to quotients (see \cite [p~264]{Abr78}) implies that it induces a smooth map $\hat{B^{\flat}}:\mathfrak{g}/G\rightarrow\mathfrak{g}^{*}/G$, $\hat{B^{\flat}}[X] = [\alpha]:=[B^{\flat}(X)]$, where $[X]$ is adjoint orbit through $X$ and $[\alpha]:=[B^{\flat}(X)]$ the corresponding coadjoint orbit through $B^{\flat}(X)=\alpha$. This gives the following diagram
\[
\xymatrixcolsep{4pc}\xymatrixrowsep{4pc}
\xymatrix{
G\times\mathfrak{g}\ar[d]_{Ad} \ar[r]^-{u} &G\times\mathfrak{g}^{*}\ar[d]^{Ad^{*}}\\
\mathfrak{g} \ar[d]_{\pi_{\mathfrak{g}}} \ar[r]^-{B^{\flat}} &\mathfrak{g}^{*}\ar[d]^{\pi_{\mathfrak{g}^{*}}}\\
\mathfrak{g}/G \ar[r]^{\hat{B^{\flat}}} &\mathfrak{g}^{*}/G}
\]

\begin{theorem}
Let $G$ be a compact, connected semisimple Lie group. Let $\mathfrak{g}$ be its Lie algebra and $\mathfrak{g}^{*}$ the dual of $\mathfrak{g}$. Let $B^{\flat}$ be as in Theorem 5.0.1 and let $\hat{B^{\flat}}:\mathfrak{g}/G\rightarrow \mathfrak{g}^{*}/G$ be the map induced by passage to quotients as described above between adjoint and coadjoint orbit spaces. Then the map $\hat{B^{\flat}}$ is a local symplectomorphism.
\end{theorem}

\proof 
The map $\hat{B^{\flat}}$ is well defined since if $\hat{B^{\flat}}([X])= [B^{\flat}(X)]$ and  $\hat{B^{\flat}}([X])=[B^{\flat}(Y)]$, then $X$ and $Y$ belong to the same orbit $[X]$ so that there is some $g\in G$ such that $Y = gXg^{-1}$. Let $\alpha = B^{\flat}(X)$ and $\beta = B^{\flat}(Y)$. Then $\beta = B^{\flat}(Y) = B^{\flat}(gXg^{-1})=gB^{\flat}(X)g^{-1}=g\alpha g^{-1}$. This shows that $\alpha$ and $\beta$ belong to the same orbit. Therefore, $[B^{\flat}(X)]= [B^{\flat}(Y)]$ so that $\hat{B^{\flat}}$ is well defined.\\

 To show that $\hat{B^{\flat}}$ is injective we first have to show that the following diagram commutes.

 \[
\xymatrixcolsep{4pc}\xymatrixrowsep{4pc}
\xymatrix{
\mathfrak{g} \ar[d]_{\pi_{\mathfrak{g}}} \ar[r]^-{B^{\flat}} & \mathfrak{g}^{*} \ar[d]^{\pi_{\mathfrak{g}^{*}}}\\
\mathfrak{g}/G \ar[r]^{\hat{B^{\flat}}} & \mathfrak{g}^{*}/G}
\]

The commuting of this diagram is now a consequence of the fact that $B^{\flat}$ is both an isomorphism and is equivariant with respect to the adjoint action and the coadjoint action. That is, $B^{\flat}\circ Ad_{g}(X) = Ad_{g}^{*}\circ B^{\flat}(X)$ for all $X\in\mathfrak{g}$ and for all $g\in G$. If we fix $X\in\mathfrak{g}$ and let $g$ run through all the elements of $G$ then on the left we get all the elements in the orbit through $X$ while on the right we get all the elements in the orbit through $B^{\flat}(X) = \alpha$. Consequently, we must have $\hat{B^{\flat}}\circ\pi_{\mathfrak{g}}(X)=\pi_{\mathfrak{g}^{*}}\circ B^{\flat}(X)$ for all $X\in\mathfrak{g}$.

We can now show that $\hat{B^{\flat}}$ is injective. The commuting of the above diagram says that $\hat{B^{\flat}}\circ\pi_{\mathfrak{g}} = \pi_{\mathfrak{g}^{*}}\circ B^{\flat}$. Suppose $\hat{B^{\flat}}([X]) = \hat{B^{\flat}}([Y])$, then $\pi_{\mathfrak{g}^{*}}\circ B^{\flat}(X) = \pi_{\mathfrak{g}^{*}}\circ B^{\flat}(Y)$. This implies that there is a $g\in G$ such that $B^{\flat}(Y)=gB^{\flat}(X)g^{-1}$. Then for all $Z\in\mathfrak{g}$ we have $B^{\flat}(Y)Z = gB^{\flat}(X)Z)g^{-1}\Rightarrow B(Y,Z) = gB(X,Z)g^{-1}\Rightarrow B(Y,Z)=B(X,Z)\Rightarrow Y=X$ so that $[X] = [Y]$ and $\hat{B^{\flat}}$ is injective. From the relation $\hat{B^{\flat}}\circ\pi_{\mathfrak{g}} = \pi_{\mathfrak{g}^{*}}\circ B^{\flat}$, the right hand side is a composition of smooth map and on the left $\pi_{\mathfrak{g}}$ is smooth, this then implies that $\hat{B^{\flat}}$ must be a smooth map.\\
To show that $\hat{B^{\flat}}$ is a surjective map consider the following commutative diagram:

\[
\xymatrixcolsep{4pc}\xymatrixrowsep{4pc}
\xymatrix{
\mathfrak{g} \ar[d]_{\pi_{\mathfrak{g}}} \ar[rd]^{\varphi} \ar[r]^-{B^{\flat}} & \mathfrak{g}^{*} \ar[d]^{\pi_{\mathfrak{g}^{*}}}\\
\mathfrak{g}/G \ar[r]^{\hat{B^{\flat}}} & \mathfrak{g}^{*}/G}
\]

We have $\varphi = \pi_{\mathfrak{g}^{*}}\circ B^{\flat}$. But the right hand side is surjective since $B^{\flat}$ is an isomorphism hence bijective and $\pi_{\mathfrak{g}^{*}}$ is the projection which is surjective, this shows that $\varphi:\mathfrak{g}\rightarrow\mathfrak{g}^{*}/G$, $X\mapsto [B^{\flat}(X)]$ is surjective. But $\hat{B^{\flat}}$ is the factorization of $\varphi$ through $\mathfrak{g}/G$,(see also \cite [pp~15-16]{Ton64}), that is, $\varphi = \hat{B^{\flat}}\circ \pi_{\mathfrak{g}}$, therefore, for any $ [B^{\flat}(X)]\in \mathfrak{g}^{*}/G$ there is $X\in\mathfrak{g}$ such that $\varphi(X)=[B^{\flat}(X)]$. This gives $\varphi(X)=\hat{B^{\flat}}(\pi_{\mathfrak{g}}(X)) = \hat{B^{\flat}}([X]) = [B^{\flat}(X)]$. Thus for each $[B^{\flat}(X)]\in\mathfrak{g}^{*}/G$ there is $[X]\in\mathfrak{g}/G$ such that $\hat{B^{\flat}}([X]) = [B^{\flat}(X)]$ which shows that $\hat{B^{\flat}}$ is bijective so that its inverse $(\hat{B^{\flat}})^{-1}$ exists. We must show that the inverse is smooth. But now $(\hat{B^{\flat}})^{-1}\circ\pi_{\mathfrak{g}^{*}}\circ B^{\flat} = \pi_{\mathfrak{g}}$ and since $\pi_{\mathfrak{g}}$ is smooth and the other two maps on the left are smooth, this forces $(\hat{B^{\flat}})^{-1}$ to be smooth. Therefore, $\hat{B^{\flat}}$ is a diffeomorphism. We shall now write $O_{X}$ for the orbit $[X]$ and $O_{B^{\flat}(X)}$ for the orbit $[B^{\flat}(X)]$.\\

Let $O_{X}$ be the adjoint orbit through $X\in\mathfrak{g}$. Define a set map on $O_{X}$ as follows: Since each element in $O_{X}$ is of the form $gX$ for some $g\in G$, for any two points $y=hX$ and $z=gX$ in $O_{X}$ let
\[
f_{X}:O_{X}\rightarrow O_{X},\quad y\mapsto z;\quad f_{X}(y)=(gh^{-1})y=z.
\]

Then $f_{X}$ maps all points of $O_{X}$ into points of $O_{X}$. Since $G$ is a group and $gh^{-1}$ is smooth for all $g,h\in G$, the map $f_{X}$ is smooth with smooth inverse $f_{X}^{-1}=hg^{-1}$.\\

In a similar way define a set map $k_{\alpha}$ on the coadjoint orbit  $O_{B^{\flat}(X)}= O_{\alpha}$ corresponding to the adjoint orbit $O_{X}$. That is,
\[
k_{\alpha}:O_{\alpha}\rightarrow O_{\alpha},\quad \beta\mapsto\gamma;\quad k_{\alpha}(\beta)=(rs^{-1})\beta=\gamma,
\]

where $\alpha = B^{\flat}(X), \beta = s\alpha, \gamma = r\alpha$ and $r,s\in G$. Let $\hat{B^{\flat}}_{X}$ be the restriction of $\hat{B^{\flat}}$ to a small neighborhood of the point $O_{X}$. Then
\[
\begin{array}{ccc}
k_{\alpha}\circ \hat{B_{X}^{\flat}}\circ f_{X}^{-1}:O_{X}\rightarrow O_{B^{\flat}(X)} = O_{\alpha}\hspace{2cm}(1)
\end{array}
\]
maps points of $O_{X}$ into points of $O_{B^{\flat}(X)}=O_{\alpha}$ and it is smooth since it is a composition of smooth maps. It is known that the coadjoint orbit is symplectic. Let $\hat{\omega}$ be the Kirillov-Kostant-Souriau form on the coadjoint orbit $O_{B^{\flat}(X)}=O_{\alpha}$ which is known to be symplectic. Then for all $Y,Z\in \mathfrak{g}$ and $r,s\in G$ we have:
\[
\begin{array}{cll}

k_{\alpha}^{*}\hat{\omega}(Y,Z) &=& \hat{\omega}(k_{\alpha{*}}Y,k_{\alpha{*}}Z)\\
        						&=& \hat{\omega}\left((rs^{-1})_{*}Y,(rs^{-1})_{*}Z\right)\\
        						&=& \hat{\omega}\left(r_{*}(s_{*}^{-1}Y),r_{*}(s_{*}^{-1}Z)\right)\\
        						&=& \hat{\omega}(r_{*}Y,r_{*}Z)\\
        						&=& \hat{\omega}(Y,Z)
\end{array}
\]
since $Y,Z\in\mathfrak{g}$ are left invariant. Thus $k_{\alpha}^{*}\hat{\omega} = \hat{\omega}$. By similar calculations, for any 2-form $\hat{\Omega}$ on the adjoint orbit $O_{X}$ we must have $f_{X}^{*}\hat{\Omega} = \hat{\Omega}$.

Consider now the pull back of the form $\hat{\omega}$ by the map in (1), $\left(k_{\alpha}\circ \hat{B_{X}^{\flat}}\circ f_{X}^{-1}\right)^{*}\hat{\omega}$.  We have
\[
\begin{array}{cll}
\left(k_{\alpha}\circ \hat{B_{X}^{\flat}}\circ f_{X}^{-1}\right)^{*}\hat{\omega} &=& (f_{X}^{-1})^{*}\circ (\hat{B_{X}^{\flat}})^{*}\circ k_{\alpha}^{*}\hat{\omega}\\
    &=& (f_{X}^{-1})^{*}\circ (\hat{B_{X}^{\flat}})^{*}\hat{\omega}
\end{array}
\]

But $\hat{B_{X}^{\flat}}$ is a smooth map so that it pulls back a 2-form into a 2-form. Thus $(\hat{B_{X}^{\flat}})^{*}\hat{\omega}$ is a 2-form. We now check if the 2-form $(\hat{B_{X}^{\flat}})^{*}\hat{\omega}$ is symplectic, that is, if it is closed and nondegenerate. Since a pull back commutes with exterior derivative we have $d\hat{B_{X}^{\flat{*}}}\hat{\omega} = (\hat{B_{X}^{\flat}})^{*}d\hat{\omega} = 0$ since $\hat{\omega}$ is closed. Thus the 2-form $(\hat{B_{X}^{\flat}})^{*}\hat{\omega}$ is closed. For non degeneracy, if $(\hat{B_{X}^{\flat}})^{*}\hat{\omega}(Y,Z)=0$ for all $Z\in\mathfrak{g}$ then $\hat{\omega}(d\hat{B^{\flat}}_{X}(Y),d\hat{B^{\flat}}_{X}(Z))=0$ for all $Z\in\mathfrak{g}$. Since $\hat{\omega}$ is symplectic, $\hat{\omega}(d\hat{B^{\flat}}_{X}(Y),d\hat{B^{\flat}}_{X}(Z))=0$ for all $Z\in\mathfrak{g}$ implies that $d\hat{B^{\flat}}_{X}(Y) = 0$. But $d\hat{B^{\flat}}$ is a linear isomorphism so that $d\hat{B^{\flat}}_{X}(Y)=0\Rightarrow Y\in\ker{d\hat{B^{\flat}}} = \lbrace 0\rbrace$ which gives $Y=0$. Thus $(\hat{B_{X}^{\flat}})^{*}\hat{\omega}(Y,Z)=0$ for all $Z\in\mathfrak{g}$ implies that $Y=0$ and $(\hat{B^{\flat}})_{X}^{*}\hat{\omega}$ is nondegenerate. This proves that $\hat{B^{\flat}}$ is a symplectic map orbitwise. So $\hat{B^{\flat}}$ can be used to pull back a symplectic form on a coadjoint orbit space to a symplectic form on an adjoint orbit space. Since the action is transitive by assumption, the orbit spaces reduce to only one each. In this case, we have proved that they are symplectomorphic spaces. More details will appear elsewhere.

\section{Acknowledgements}
Augustin Batubenge is grateful to Professor Fran\c{c}ois Lalonde for his financial support and for hosting him as an invited researcher in the Canada chair of mathematics during the time of writing this paper at the University of Montr\'eal from 2017 to 2019.\\
Wallace Haziyu acknowledges the financial support from the International Science Program, ISP, through East African Universities Mathematics Project, EAUMP and more particularly to Professor Lief Abrahamson for his significant input in funding his research.

{\bf Authors}
\begin{itemize}
\item Augustin Tshidibi Batubenge\\
	Department of Mathematics and Statistics
	Universit\'e de Montr\'eal and University of Zambia\\
	email: a.batubenge@gmail.com\\
\item Wallace Mulenga Haziyu\\
	Department of Mathematics and Statistics\\
	University of Zambia\\
	P.O. Box 32379 Lusaka, Zambia\\
	email: whaziyu@unza.zm
\end{itemize}

\end{document}